\baselineskip=14pt
\overfullrule=0mm

\hsize=118mm 
\hoffset=20mm
\vsize=215mm
\voffset=15mm

\font\tenbb=msbm10
\font\sevenbb=msbm7
\font\fivebb=msbm5
\newfam\bbfam
\textfont\bbfam=\tenbb \scriptfont\bbfam=\sevenbb
\scriptscriptfont\bbfam=\fivebb
\def\bb{\fam\bbfam}

\def\Qb{{\bb Q}}
\def\Rb{{\bb R}}
\def\Zb{{\bb Z}}

\def\D{\Delta}
\def\g{\gamma}
\def\Lb{\wedge}
\def\om{\omega}
\def\s{\sigma}
\def\Si{\Sigma}
\def\vp{\varphi}
\def\z{\zeta}

\def\bu{\cdot}
\def\ify{\infty}
\def\ot{\otimes}
\def\ov{\overline}
\def\part{\partial}
\def\sbs{\subset}
\def\wh{\widehat}
\def\lb{\lambda}

\def\ra{\rightarrow}

\def\card{\mathop{\rm card}\nolimits}
\def\Gal{\mathop{\rm Gal}\nolimits}
\def\Im{\mathop{\rm Im}\nolimits}
\def\Inf{\mathop{\rm Inf}\nolimits}
\def\Spec{\mathop{\rm Spec}\nolimits}
\def\Sup{\mathop{\rm Sup}\nolimits}
\def\tors{\mathop{\rm tors}\nolimits}
\def\Ker{\mathop{\rm Ker}\nolimits}

\def\build#1_#2^#3{\mathrel{
\mathop{\kern 0pt#1}\limits_{#2}^{#3}}}

\def\Cc{{\cal C}}
\def\Sc{{\cal S}}

\def\limproj{\mathop{\oalign{lim\cr
\hidewidth$\longleftarrow$\hidewidth\cr}}}


\vglue 1cm

\centerline{\bf Perfect forms and the Vandiver conjecture }

\bigskip

\centerline{C. SOUL\'E}

\vglue 1cm

\hfill To J.Neukirch

\vglue 1cm

Let $p$ be an odd prime, $C$ the $p$-Sylow subgroup of 
the class group of $\Qb (\root p\of{1})$, 
and $C^+$ the subgroup of $C$ fixed by complex conjugation. The Vandiver 
conjecture is the statement that $C^+ = 0$. It has been verified when $p$ is 
less than four million [B-C-E-M].

\smallskip

For any natural integer $i \leq p-2$, let $C^{(i)}$ be the subgroup of $C$ where the 
Galois group of $\Qb (\root p\of{1})$ over $\Qb$ acts by the $i$-th power
of the 
Teichm\"uller character. Vandiver's conjecture says that $C^{(i)}$ vanishes for 
all even values of $i$. Kurihara proved in [K] that 
$C^{(p-3)} = 0$. His proof uses the existence of a surjective map
$$
K_{2n-2} (\Zb) \ra C^{(p-n)} \build \ot_{\Zb}^{} \Zb /p
$$
and the fact that $K_4 (\Zb)$ is not too big. Here $K_m (\Zb)$ is the $m$-th 
algebraic $K$-group of the ring of integers, $m>0$. Since $K_m (\Zb)$ is finite 
unless $m-1$ is positive and divisible by four [B], for a fixed value of $n$ 
there are only finitely many $p$ with $C^{(p-n)} \not= 0$. In this paper we 
show that, when $n$ is odd and positive,
$$
C^{(p-n)} = 0 \quad \hbox{when} \quad p > v(n) \, ,
$$
where
$$
\log v(n) \leq n^{224n^4} \, .
$$

This huge bound follows from a similar assertion about the torsion subgroup of 
$K_m (\Zb)$ which in turn comes from a bound on the torsion of the various 
homology groups of the special linear group $G = SL_N (\Zb)$ with integral 
coefficients, $N>1$. The homology of $G$ coincides, up to finite groups of 
small cardinality, with the homology of the quotient $Y = X/G$, where $X$ is 
the symmetric space $SO_N \backslash SL_N (\Rb)$. A finite cell decomposition 
of a compactification $Y^*$ of $Y$ is provided by the beautiful Voronoi 
reduction theory [V].

\smallskip

The crucial step in our argument is the fact, that we learnt from \break 
O.~Gabber, that a bound on the number of cells of a $CW$-complex as $Y^*$ (and 
on the number of faces of any given cell in $Y^*$) provides an explicit upper 
bound for the cardinality of the torsion subgroup of its integral homology. 
Such bounds for the Voronoi cell decomposition of $Y^*$ follow from standard 
arguments on lattices. Since both steps are exponential in $N$, we cannot get 
better than a double exponential estimate for $v(n)$.

\smallskip

In Section 1 we discuss Voronoi's reduction theory. In Section 2 we bound the 
order of the torsion in the homology of $SL_N (\Zb)$ (Theorem 1) and in the 
$K$-theory of $\Zb$ (Theorem 2). In Section 3 we apply these bounds to the 
Vandiver conjecture (Theorem 4) and to groups related to it by the Iwasawa 
theory (Theorem 3 and Theorem 5).

\smallskip

It would be interesting to generalize the results above by replacing $\Zb$ by the 
ring of integers in an arbitrary number field.

\smallskip

I thank O.~Gabber for telling me about Proposition 3, which was the starting 
point of this work, as well as J.~Buhler, O.~Gabber,
J.~Martinet, J-P.~Serre, 
L.~Washington and D.~Zagier for helpul comments.

\vglue 1cm

\noindent {\bf 1. Perfect forms}

\medskip

\noindent {\bf 1.1.} Let $N \geq 2$ be an integer, let $V_N$ be the 
vector space of $N$ by $N$ real symmetric matrices, and let $P_N \sbs 
V_N$ be the open cone of all positive definite symmetric matrices. The 
group $\Rb_+$ of positive real numbers acts on $P_N$ by multiplication, and 
the quotient $X_N = P_N / \Rb_+$ is the symmetric space of $SL_N (\Rb)$.

\smallskip

A matrix $A$ in $V_N$ is said to have {\it rational nullspace} when the 
nullspace $\Ker (A)$ of the bilinear form defined by $A$ is spanned by 
vectors in $\Qb^N \sbs \Rb^N$. We let $P_N^*$ be the subset of $V_N$ 
consisting of all nonzero positive semi-definite symmetric matrices with 
rational nullspace, $X_N^* = P_N^* / \Rb_+$ and $\pi : P_N^* \ra X_N^*$ 
 the projection map.

\smallskip

The group $GL_N (\Zb)$ acts upon $P_N^*$ on the right by the formula
$$
A \cdot g = g^t \, A \, g \, , \ g \in GL_N (\Zb) \, , \ A \in P_N^* \, ,
$$
where $g^t$ is the transpose of $g$,
and $P_N$ is invariant under this action.

\smallskip

Given $A$ in $P_N$ we let $\mu (A)$ be the minimum value of the number 
$v^t \, A \, v$ where $v$ ranges over all nonzero vectors in the lattice 
$L = \Zb^N \sbs \Rb^N$, and $m(A)$ be the (finite) set of vectors $v\in 
L - \{ 0 \}$ such that $v^t\, A \, v = \mu (A)$, i.e. the set of minimal 
vectors.

\medskip

 A form $A \in P_N$ is called {\it 
perfect} if $\mu (A) = 1$ and, for any $B$ in $P_N$ such that $\mu (B) = 
1$, the equality of sets $m(B) = m(A)$ implies $B=A$.

\bigskip

\noindent {\bf 1.2.} It was shown by Voronoi [V] p.110, Thm.,
 that, up to conjugation by 
$SL_N (\Zb)$, there are only finitely many perfect forms. To give an upper 
bound on the number of classes of perfect forms,
 it is clearly enough to bound the size 
of their minimal vectors. Let us denote by $\g = \g (N)$ the {\it Hermite 
constant} i.e. the supremum of the quantity $\mu (A) \, \det (A)^{-1/N}$ 
when $A$ runs over $P_N$. We also denote by $s(N)$ the maximum of $\card 
(m(A))/2$ over all $A \in P_N$ (half the {\it kissing number}
 of the corresponding 
packing of spheres, [C-S]~I), where $\card (S)$ 
denotes the cardinal of a finite 
set $S$.

\medskip

\noindent {\bf Proposition 1.} {\it Let $A \in P_N$ be such that $\mu (A) = 
1$ and $m(A)$ spans the real vector space $\Rb^N$. Then there exists $g \in 
SL_N (\Zb)$ such that any vector $v \in m(A \cdot g)$ has 
coordinates $x_i$ such 
that
$$
\vert x_i \vert \leq A(N) \, , \quad 1 \leq i \leq N \, , \leqno (1)
$$
with
$$
A(N) = N^{(N-1)} \, \g^{N/2} \, . \leqno (2)
$$
}

\medskip

\noindent {\it Proof of Proposition 1.} (See [M], III~5.2.) Let
$$
h(v,w) = \ v^t \, A \, w \, , \quad v,w \in \Rb^N \, ,
$$
be the quadratic form defined by $A$ and $h(v) = h(v,v)$. By [Z], Lemma 
1.7, since $m(A)$ spans $\Rb^N$ and $\mu (A) = 1$ we can find a basis of 
$\Zb^N$ made of vectors $(e_i)$ such that
$$
h(e_i) \leq N^2 \, , \quad 1 \leq i \leq N \, . \leqno (3)
$$
Up to conjugation by $SL_N (\Zb)$ and, eventually, the replacement of $e_1$ 
by $-e_1$, we may assume that $(e_i)$ is the standard basis 
of $L=\Zb^N$. Let
$$
v = \sum_{i=1}^N x_i \, e_i
$$
be any vector in $\Rb^N$. For a any integer $i$, $1 \leq i \leq N$, denote by 
$A_i$ the matrix of scalar products
$$
A_i = (h (v_k , v_{\ell}))
$$
where $v_k = e_k$ when $k \not= i$ and $v_i = v$. Clearly
$$
\vert x_i \vert^2 = \det (A_i) / \det (A) \, . \leqno (4)
$$
When $v$ is in $m(A)$, Hadamard's inequality, 
(3) and the fact that $h(v)=1$ imply that
$$
\det (A_i) \leq N^{2(N-1)} \, . \leqno (5)
$$
The definition of $\g$ implies
$$
\det (A)^{-1} \leq \g^N \, . \leqno (6)
$$
The proposition follows from (4), (5) and (6).

\hfill q.e.d.
\bigskip

\noindent {\bf 1.3.}
Any vector $v$ in $L - \{ 0 \}$ determines a form $\wh v = v \, v^t$ in 
$P_N^*$. Given any finite subset $B \sbs L - \{ 0 \}$, the {\it convex 
hull} of $B$ is the image by $\pi$ of the subset 
$$
\left\{ \sum_j \lb_j \, \wh{v_j} \, : \, v_j \in B \, , \, \lb_j \geq 0 
\right\}
$$
of $P_N^*$. When $A$ is a perfect form, we let $\s (A) \sbs X_N^*$ be the 
convex hull of its set $m(A)$ of minimal vectors.

\smallskip

According to Voronoi [V], \S 8-15, see also [As],
the cells $\s (A)$ and their intersections, when $A$ 
runs over all perfect forms, define a cell decomposition of $X_N^*$, 
invariant under $GL_N (\Zb)$. We equip $X_N^*$ with the corresponding 
$CW$-topology. Since there are only finitely many perfect forms modulo the 
action of $SL_N (\Zb)$ ,  the quotient space $Y_N^* = X_N^* / 
SL_N (\Zb)$ is a finite $CW$-complex. 

\medskip

Let now
$$
B(N) = (2 A(N) + 1)^N / 2 \, , \leqno (7)
$$
 where $A(N)$ is defined by (2), and let $k \geq 0$ be an integer.

\medskip

\noindent {\bf Proposition 2.}

\smallskip

\item{i)} {\it The number of equivalence classes of $k$-dimensional cells 
in the Voronoi decomposition of $X_N^*$ is at most
$$
c (k,N) = \left( \matrix{ B(N) \cr k+1 } \right) \, . \leqno (8)
$$
}

\item{ii)} {\it Any $k$-dimensional cell has at most $f(k,N)$ faces, with 
$$
f(k,N) = \left( \matrix{ s(N) \cr k } \right) \, . \leqno (9)
$$
}

\medskip

\noindent {\it Proof.} Let $\Phi$ be the set of non zero
vectors $v = (x_i)$ in 
$\Zb^N$ satisfying (1). Clearly $\card (\Phi) \leq 2B(N)$.  
Any cell  $\s$ is the convex hull
of  a subset of $m(A)$, for some perfect form
$A$ ( [V], \S 20, [M], VII Thm. 1.12). Therefore,
if it has dimension $k$
 there exists $g \in SL_N (\Zb)$ 
such that the interior of
$\tau = \s \cdot g$ contains the interior of the convex hull of 
$(k+1)$ vectors in $\Phi$. There are at most $c (k,N)$ such cells $\tau$.

\medskip

For the same reason,
any face of $\s$ contains the interior of the convex hull of $k$ 
vectors in $m(A)$. The number of such $k$-uples 
of vectors (taken up to sign) is bounded by 
$f(k,N)$.

\hfill q.e.d.

\bigskip

\noindent {\bf 2. Homology of $SL_N (\Zb)$}

\medskip

\noindent {\bf 2.1.} Let $a>0$, $b>0$ be integers and 
$$
\vp : \Zb^a \ra \Zb^b
$$
a $\Zb$-linear map. Denote by $S$ the image of $\vp$, by $Q$ the cokernel 
of $\vp$ and by $Q_{\tors}$ the torsion subgroup of $Q$. Denote by $\Vert 
\cdot \Vert$ the standard euclidean norm in $\Rb^b$ and by $(e_i)$, $1 \leq 
i \leq a$, the standard basis of $\Zb^a$.

\medskip

\noindent {\bf Lemma 1.} {\it Let $I \sbs \{ 1, \ldots , a \}$ be a set of 
indices such that the set of vectors $\vp (e_i)$, $i \in I$, is a basis of 
$S \build \ot_{\Zb}^{} \Rb$. Then
$$
\card (Q_{\tors}) \leq \prod_{i \in I} \Vert \vp (e_i) \Vert \, .
$$
}

\medskip

\noindent {\it Proof.} If $E = \Zb^b$ we have an exact sequence of finitely 
generated abelian groups
$$
0 \ra S \ra E \ra Q \ra 0 \, . \leqno (10)
$$
Equip $E \build \ot_{\Zb}^{} \Rb = \Rb^b$ with the standard scalar product, 
and both $S \ot \Rb$ and $Q \ot \Rb$ with the induced scalar product (so 
that the projection map $E \ot \Rb \ra Q \ot \Rb$ gives an isometry from 
the orthogonal complement of $S \ot \Rb$ onto $Q \ot \Rb$). Call $\ov{S}$, 
$\ov{E}$, $\ov{Q}$ these euclidean lattices. Their arithmetic degrees ([L] 
V \S~2, [S4] VIII 1.4) satisfy the relation
$$
\wh{\deg} (\ov{E}) = \wh{\deg} (\ov{S}) + \wh{\deg} (\ov{Q}) \, . \leqno 
(11)
$$
Since $\ov E$ is standard we have
$$
\wh{\deg} (\ov{E}) = 0 \, . \leqno (12)
$$
Let $P = Q / Q_{\tors}$. Since $P$ is torsion free and spanned by 
vectors of length at most one, we get
$$
\wh{\deg} (\ov{P}) = -\log {\rm covolume} \, (\ov{P}) \geq 0 \, .
$$
Therefore
$$
\wh{\deg} (\ov{Q}) = \wh{\deg} (\ov{P}) + \wh{\deg} (Q_{\tors}) \geq \log 
\card (Q_{\tors}) \, . \leqno (13)
$$ 
Furthermore, the exterior product $\build \Lb_{ i\in I}^{} \vp (e_i)$ 
(taken in any order) is a non zero element of the maximal exterior power 
$\det (S)$ of the lattice $S$. Therefore, by the Hadamard inequality,
$$
\wh{\deg} (\ov{S}) = \wh{\deg} (\det (\ov{S})) \geq - \log \Vert \build 
\Lb_{ i\in I}^{} \vp (e_i) \Vert \geq - \sum_{i \in I} \log \Vert \vp (e_i) 
\Vert \, . \leqno (14)
$$
>From (11), (12), (13) and (14) we conclude that
$$
\log \card (Q_{\tors}) \leq \sum_{i \in I} \log \Vert \vp (e_i) \Vert \, .
$$

\hfill q.e.d.

\bigskip

\noindent {\bf 2.2.} Let $(C_{\bu},\part)$ be a chain complex of free 
finitely generated $\Zb$-modules $C_k$, $k \geq 0$. Let $\Si_k$ be a basis 
of $C_k$ and, for any $\s \in \Si_{k+1}$, let us write
$$
\part (c) = \sum_{c' \in \Si_k} n_{\s \s'} \, \s' \, .
$$
Define
$$
a(k) = \Inf (\card (\Si_{k+1}) , \card (\Si_k )) \leqno (15)
$$
and
$$
b(k) = \Sup \left( \build \Sup_{\s \in \Si_{k+1}}^{} \left( \sum_{\s' \in 
\Si_k} n_{\s \s'}^2 \right)^{1/2} , 1 \right) \, . \leqno (16)
$$

\medskip

\noindent {\bf Proposition 3.} (O. Gabber): {\it For any $k \geq 0$,
$$
\card H_k (C_{\bu})_{\tors} \leq b(k)^{a(k)} \, .
$$
}

\medskip

\noindent {\it Proof.} The homology group
$$
H_k (C_{\bu}) = {{\rm Ker} (\part) \over \Im (\part)}
$$
is contained in the cokernel $Q$ of the map
$$
\part : C_{k+1} \ra C_k \, .
$$
For any $\s \in \Si_{k+1}$, the vector $\part (\s)$ has coordinates $(n_{\s 
\s'})$ in the free $\Zb$-module $C_k = \Zb^{\Si_k}$. Therefore
$$
\Vert \part (\s) \Vert = \left( \sum_{\s'} n_{\s \s'}^2 \right)^{1/2} \leq b 
(k) \, .
$$
The rank $t$ of the image of $d$ is at most $a(k)$. From Lemma 1 we conclude 
that
$$
\card H_k (C_{\bu})_{\tors} \leq \card Q_{\tors} \leq b(k)^t \leq 
b(k)^{a(k)} \, .
$$
\hfill q.e.d.

\bigskip

\noindent {\bf 2.3.} For any integer $n > 0$ we let $\Sc_n$ be the Serre 
class of finite abelian groups which have no element of prime order $p > 
n$. Given any finite abelian group $A$, let $B \sbs A$ be the maximal 
subgroup of $A$ lying in $\Sc_n$ and
$$
\card_n (A) = \card (A/B) \, .
$$
This quantity depends only on the class of $A$ modulo $\Sc_n$.

\bigskip

\noindent {\bf 2.4.} Let $N$ and $k$ be positive integers, with $k \leq 
N(N-1) / 2$. Let
$$
\ell = (N(N+1)/2) - 1 - k \, ,
$$
and
$$
h(k,N) = f(\ell , N)^{c(\ell , N)} \, , \leqno (17)
$$
where $f$ and $c$ are defined by (8) and (9).

\medskip

\noindent {\bf Theorem 1.} {\it The torsion subgroup in the homology of 
$SL_N (\Zb)$ is bounded as follows:
$$
\card_{N+1} H_k (SL_N (\Zb) , \Zb)_{\tors} \leq h (k,N) \, .
$$
}

\medskip

\noindent {\it Proof.} We may assume $N \geq 3$ since the homology of $SL_2 
(\Zb)$ is well-known (and we do not need it).
 Let $G = SL_N (\Zb)$, $X^* = X_N^*$
 and $ \part X^* = \part X_N^*$.
Denote by $C_{\bu} (X^* , \part X^* )$  
the chain complex of 
the pair of $CW$-complexes $(X^* , \part X^* )$ 
with integral coefficients, 
by $C_{\bu} (X^* , \part X^* )_G$  its coinvariants 
under the action of $G$, and by
$(C_{\bu} , \part )= C_{\bu} (X^* , \part X^* )_G/(\rm{torsion})$ 
the quotient of this complex by its torsion sub-complex. 
For all $k \geq 0$, let $\Si_k'$ be a set of representatives, 
modulo $G$, of the $k$-dimensional cells in 
the Voronoi cell decomposition 
of $X^*$ which are not contained in $\part X^*$. One has
$$
C_k (X^* , \part X^* )_G
\sim \bigoplus_{\s \in \Si_k'} H_0 (G_{\s} , \Zb_{\s}),
$$
where $G_{\s}$ is the stabilizer of $\s$ and $ \Zb_{\s}$
its orientation module. Denote by $\Si_k \subset \Si_k'$
the set of cells $\s$ such that no element in $G_{\s}$
changes the orientation of $\s$. When $\s$ lies in
 $\Si_k $ the group $H_0 (G_{\s} , \Zb_{\s})$ is isomorphic
to $\Zb$, and it is killed by two if $\s$ lies in 
$\Si_k' - \Si_k$. Therefore $C_k$ is a free module 
with basis $\Si_k$ over $\Zb$.

 Given any $\s \in 
\Si_{k+1}$ we have
$$
\part (\s) = \sum_{\s' \in \Si_k} n_{\s \s'} \, e_{\s'},
$$
where $\vert n_{\s \s'} \vert$ is at most the number of faces $\tau$ of 
$\s$ in $X^*$ which are equivalent to $\s'$. Proposition 2 implies that, 
for any $\s \in \Si_{k+1}$,
$$
\left( \sum_{\s' \in \Si_k} n_{\s \s'}^2 \right)^{1/2} \leq \sum_{\s' \in 
\Si_k} \vert n_{\s \s'} \vert \leq f(k+1,N) \leqno (18)
$$
and
$$
\card (\Si_k) \leq c (k,N) \, . \leqno (19)
$$
>From (18), (19) and Proposition 3 we conclude that
$$
\card H_k (C_{\bu})_{\tors} \leq f(k+1,N)^{c(k+1,N)} \, . \leqno (20)
$$
Since $C_{\bu}$ is torsion free, by the universal coefficient
theorem, we have
$$
\card H^{k+1} (C_{\bu} , \part)_{\tors}
 = \card H_k (C_{\bu} , \part)_{\tors} \, . 
\leqno (21)
$$

On the other hand, given any element $g$ in the group $G$ of 
order a prime
$p $,  the roots of the characteristic
polynomial of $g$ are $p$-th roots of unity,
therefore $p \leq N+1$.
 Furthermore, the stabilizer of any cell of $X^*$
which is not contained in $\part X^*$ is finite. 
Therefore the cohomology group $H^k (C_{\bu} , \part)$
 coincides modulo 
$\Sc_{N+1}$ with the equivariant
 cohomology group $H_G^k (X^* , \part X^*)$ 
([Br], VII.7). From [S5] Proposition 1, we know that
$$
H^m (X^* , \part X^*) = \left\{ \matrix{
St^* \ {\rm if} \ m = N-1 \, , \hfill \cr
0 \ {\rm otherwise}, \hfill \cr
} \right.
$$
where $St^*$ is the top cohomology of the Borel-Tits building of $SL_N$ over 
$\Qb$, i.e. the $\Zb$-dual of the Steinberg module. Therefore 
( by the cohomological analog of [Br] VII 
(7.2))
$$
H_G^k (X^* , \part X^*) = H^{k-N+1} (G,St^*) \, . \leqno (22)
$$
The long exact sequence of Farrell, [F] Theorem 2 (b), and the Borel-Serre 
duality [B-S] tell us that
$$
H^k (G,St^*) = H_{v-k} (G,\Zb) \leqno (23)
$$
modulo $\Sc_{N+1}$, where $v = N (N-1)/2$ is the virtual cohomological 
dimension of $G$. Indeed, since   no element in the group $G$
 has order a prime $p > N+1$, 
the Farrell homology of $G$
lies in $\Sc_{N+1}$, as follows from
the homological analog of the spectral sequence
 of   [Br] X (4.1) (one can also use
 [Br] X. 3 Exercise 2).
  In particular, $H_k (G,\Zb)$ lies in $\Sc_{N+1}$ 
when $k > 
N(N-1)/2$ and  the theorem would be trivial in that case. If 
we combine (20), (21), (22) and (23) we get our assertion.

\hfill q.e.d.
\bigskip

\noindent {\bf 2.5.} For any $m \geq 1$ let $K_m (\Zb)$ be the $m$-th 
higher algebraic $K$-group of the integers, and
$$
k(m) = h (m,2m+1) \, . \leqno (24)
$$

\medskip

\noindent {\bf Theorem 2.} {\it For any $m > 0$ we have
$$
\card_{2m+2} K_m (\Zb)_{\tors} \leq k(m) \, .
$$
}

\medskip

\noindent {\it Proof.} Consider the Hurewicz map of the $H$-space $BGL 
(\Zb)^+$:
$$
H : K_m (\Zb) \ra H_m (GL (\Zb) , \Zb) \, .
$$
The kernel of $H$ lies in $\Sc_n$, where $n$ is the integral part of 
$(m+1)/2$ ([S2] Proposition 3, see also [A] Theorem 1.5). Furthermore
$$
H_m (GL(\Zb) , \Zb) = H_m (GL_N (\Zb) , \Zb)
$$
whenever $N \geq 2m+1$ (cf. for instance [Su], Corollary 8.3).
 When $N$ is odd, $GL_N (\Zb)$ is the product of $SL_N (\Zb)$ by a 
group of order two, so their homology groups coincide modulo $\Sc_2$. We 
now apply Theorem 1 and we get Theorem 2.

\hfill q.e.d.
\bigskip

\noindent {\bf 2.6.} When $m \leq 5$, $K_m (\Zb)$ has no $p$-torsion unless 
$p \leq 3$ [L-S]. Let us evaluate $k(m)$ when $m \geq 6$.

\medskip

\noindent {\bf Lemma 2.} {\it When $m \geq 6$
$$
\log k(m) \leq m^{20m^4} \, .
$$
}

\medskip

\noindent {\it Proof.} The following estimates hold:
$$
\g (N) \leq 1 + {N \over 4} \leqno (25)
$$
([M-H] II (1.5) Remark) and
$$
s(N) \leq 2^N - 1\leqno (26)
$$
([V] p.107, Lemme) (for sharper estimates,
 see [K-L] and [C-S] Chapter I, (49) and 
(50)). When $N = 2m+1$ and $k=m$, we have
$$
\ell = {N (N+1) \over 2} - m - 1 = (N^2 - 1)/2 \, .
$$
By the definitions (17) and (24) we have
$$
\log \log k(m) = \log c (\ell , N) + \log \log f(\ell , N) \, . \leqno (27)
$$
Since
$$
c(\ell , N) = \left( \matrix{ B(N) \cr \ell + 1 } \right) \leq B(N)^{\ell + 1}
$$
we get
$$
\log c (\ell , N) \leq {N^2 + 1 \over 2} \, \log (B(N)) \, . \leqno (28)
$$
Using (7), (2) and (25) we have
$$
\eqalign{
B(N) = & \ (2A(N)+1)^N / 2\, , \cr
A(N) \leq & \ N^{(N-1)} \left( 1 + {N \over 4} \right)^{N/2} \, . 
}
$$
Since
$$
\log (1+x) \leq \log (x) + {1 \over x} \, , \leqno (29)
$$
we get an upper bound of $\log B(N)$ by a finite linear combination of 
$\log (N) N^k$ and $N^k$, $-1 \leq k \leq 4$. Applying (29) again to 
replace $N$ by $m$, we get from (26), (9), (27) and (28) a similar upper 
bound for $\log \log k(m)$, namely
$$
\log \log k(m) \leq 12 m^4 \log (m) + \varepsilon(m) \leqno (30)
$$
where
$$\varepsilon(m) =
4 \log(2) m^4 + (20 \log(m)+ 14 +8 \log(2) ) m^3
     + (15 \log(m) + 22 + 7 \log(2) ) m^2 $$
$$ + (5 \log(m)+ 31/2 + 3 \log(2) ) m
     + 7 \log(m)/2 + 7 \log(2)/2  + 5.$$
Since
$$
 \varepsilon(m)\leq 8 m^4 \log (m)
$$
as soon as $m \geq 6$, Lemma 2 follows.

\hfill q.e.d.

\vglue 1cm

\noindent {\bf 3. Iwasawa theory}

\medskip

\noindent {\bf 3.1.} Let $p$ be an odd prime, $k \geq 0$ an integer and $n 
\in \Zb$. Denote by
$$
H^k (\Zb [1/p] , \Zb_p (n)) = \limproj_{\nu} H^k (\Spec (\Zb [1/p]) , \Zb / 
p^{\nu} (n))
$$
the \'etale cohomology groups 
of the scheme $\Spec (\Zb [1/p])$ with coefficients in 
the $n$-th Tate twist of the group of $p$-adic integers. When $n \not= 0$ 
these groups vanish unless $k=1$ or $2$. It was shown in [S1] and [D-F] 
that when  $m = 2n-k > 1$ and $k=1$ or $2$, there is a surjective map
$$
K_m (\Zb) \otimes \Zb_p \ra H^k (\Zb [1/p] , \Zb_p (n)) \, .
$$
(In [S1] a Chern class map is defined, the cokernel of which lies in 
$\Sc_n$. This is enough for our purpose since Theorem 2 deals only with big 
primes. In [D-F] a Chern character map is defined, which is always 
surjective.) Therefore Theorem 2 gives an upper bound for the torsion of 
these \'etale cohomology groups. However, when $k=1$ one has
$$
H^1 (\Zb [1/p] , \Zb_p (n))_{\tors}
 = H^0 (\Zb [1/p] , \Qb_p / \Zb_p (n)) = H^0 
(\Gal (\ov{\Qb} / \Qb) , \Qb_p / \Zb_p (n)) \, ,
$$
and this group is zero unless $p \leq n+1$. When $k=2$ and $n > 0$ is even, 
the main conjecture, proved by Mazur and Wiles [M-W], tells us that the 
order of $H^2 (\Zb [1/p] , \Zb_p (n))$ is the $p$-part of the numerator 
$N_n$ of $B_n / n$, where the Bernoulli numbers $B_n$
are defined by the identity 
of formal power series
$$
{x \over e^x - 1} = \sum_{n=0}^{\ify} B_n {x^n \over n!} \, .
$$
So the interesting case is when $k=2$ and $n$ is odd, i.e. when $m$ is 
divisible by 4. The group $K_m (\Zb)$ is then finite [B] and we get from 
Theorem 2 the following:

\medskip

\noindent {\bf Theorem 3.} {\it Let $n \geq 3$ be an odd integer. Then
$$
\prod_{p \geq 4n-1} \card H^2 (\Zb [1/p] , \Zb_p (n)) \leq k(2n-2) \, .
$$
}

\medskip

\noindent {\bf 3.2.} Let $\Qb (\mu_p)$ be the cyclotomic extension of $\Qb$ 
obtained by adding $p$-th roots of unity. The Galois group $\D = \Gal (\Qb 
(\mu_p) / \Qb)$ is isomorphic to $(\Zb / p)^*$, and we let
$\om : \D \ra \Zb_p^*$ be the  Teichm\"uller character. It is such
that
$$
g(\z) = \z^{\om (g)},
$$
 for any $g \in \D$ and any $p$-th root of unity $\z$.
Let $C$ be the $p$-Sylow subgroup of the class group of $\Qb 
(\mu_p)$. For any $j \in \Zb$, denote by $C^{(j)}$ the set of elements $x 
\in C$ such that
$$
g(x) = \om (g)^j \, x
$$
for any $g \in \D$. Since $\D$ has order prime to $p$, the subgroup $C^+$ 
of $C$ fixed by the complex conjugation is the direct sum of those groups 
$C^{(j)}$ such that $j$ is even and $0 \leq j \leq p-3$. One has [K]
$$
C^{(0)} = C^{(p-3)} = 0 \, .
$$
A short computation ([K], Lemma 1.1) gives that
$$
H^2 (\Zb [1/p] , \Zb_p (n)) \build \ot_{\Zb}^{} \Zb /p 
= C^{(p-n)} \build \ot_{\Zb}^{} \Zb /p\, .
$$
So Theorem 3 implies

\medskip

\noindent {\bf Theorem 4.} {\it Let $n > 1$ be an odd integer. If $p > v(n) 
= k(2n-2)$, one has
$$
C^{(p-n)} = 0 \, .
$$
}

\medskip

\noindent From the proof of Lemma 2 one gets
$$
\log \log v(n) \leq 192 n^4 \log (n) + \varepsilon(n)
$$
with
$$
\varepsilon(n)
 \leq 32 n^4 \log (n) \, ,
$$
so that $\log v(n) \leq n^{224n^4}$, if $n \geq 5$.

\bigskip

\noindent {\bf 3.3.} Let $L(s,\om^k)$ be the Dirichlet $L$-function of the 
character $\om^k$, $k \in \Zb$. Denote by $E$ the group of units in $\Qb 
(\mu_p)$ and by $E/p$ its quotient by $p$-th powers.

\medskip

\noindent {\bf Theorem 4.} {\it Let $n \geq 5$ be an odd integer, $p > 
v(n)$, and $m \in \Zb$ an integer congruent to $n$ modulo $(p-1)$.}

\item{i)} {\it The group $C^{(n)}$ is isomorphic to $\Zb_p / L (0,\om^{-n}) 
\  \Zb_p$.}

\item{ii)} {\it The cyclic group $(E/p)^{(p-n)}$ is spanned by the 
cyclotomic unit $E_{p-n}$ defined in {\rm [W]} 8.3, p.~155.}

\item{iii)} {\it The group $H^2 (\Zb [1/p] , \Zb_p (m))$ vanishes and $H^2 
(\Zb [1/p] , \Zb_p (p-m))$ is a cyclic $\Zb_p$-module.}

\item{iv)} {\it When $p > m$ the group $H^1 (\Zb [1/p] , \Zb_p (p-m))$ 
vanishes. When $m>1$ the group
$H^1 (\Zb [1/p] ,$ $\Zb_p (m))$ is a free cyclic
$\Zb_p$-module spanned 
by a cyclotomic element.}

\medskip

\noindent {\it Proof.} The $p$-rank of $C^{(n)}$ is less or equal to 1 by 
Theorem 3 and [W], Theorem 10.9. Its order is known by [M-W]. This proves 
i). 

\smallskip

It is shown in [M-W], 1.10, Theorem 1, that, if
$\Cc$ is the group of cyclotomic units, the abelian group
$((E/\Cc)\ot \Zb_p)^{(p-n)}$ has the same cardinal as $C^{(p-n)}$. 
 Therefore ii) is also a 
consequence of Theorem 3. 

Assertion iii) follows from the equality
$$
H^2 (\Zb [1/p] , \Zb_p (j)) \ot \Zb / p = C^{(p-j)} \build \ot_{\Zb}^{} \Zb /p
\, \, \, \, , j \in \Zb, $$
 that we mentionned already. 
 
 When $p > m$ the 
cyclic group $H^2 (\Zb [1/p] , \Zb_p (p-m))$ is finite, isomorphic to $\Zb_p / 
N_{p-m} \Zb_p$. 
The exact sequences
$$
0 \ra H^1 (\Zb [1/p] , \Zb_p (j)) \ot \Zb / p \ra  \ H^1 (\Zb [1/p] , \Zb / p 
(j)) 
\ra  \ H^2 (\Zb [1/p] , \Zb_p (j)) [p] \ra 0
$$
and
$$
0 \ra (E/p)^{(p-j)} \ra H^1 (\Zb [1/p] , \Zb / p (j)) \ra C [p]^{(p-j)} \ra 0 
\, ,
$$
valid when $j \not\equiv 1$ modulo $(p-1)$, and the previous discussion on 
$H^2$ prove that 
$H^1 (\Zb [1/p] , \Zb_p (m))$ is always 
cyclic over $\Zb_p$ and that, 
 if $p > m$, the group $H^1 (\Zb [1/p] , \Zb_p (p-m))$ vanishes.
 When $m>1$,  $H^1 (\Zb [1/p] , \Zb_p (m))$ is spanned by a 
cyclotomic element since the index of cyclotomic elements [S3] is the 
cardinality of $H^2 (\Zb [1/p] , \Zb_p (m))$ ([B-K], (6.8) and (6.9)).

\hfill q.e.d.

\bigskip

\noindent {\bf 3.4.} By the argument of [W], Theorem 8.14, together with [M-W], 
1.10 Theorem 1, $C^{(p-n)}$ is different from zero if and only 
if the cyclotomic unit 
$E_{p-n}$ of 
[W], 8.3, is a $p$-th power. Let us assume that this happens with probability 
$1/p$ (compare [W] loc.cit., Remark), and fix an odd integer $n$. The 
probability that there exists a prime $p \leq x$, with $p \geq 37$ (the first 
irregular prime), such that $C^{(p-n)} \not= 0$ is thus bounded above by
$$
\sum_{37 \leq p \leq x} {1 \over p} \sim \log \log (x) - 2.56 \, . 
$$
When $x = 4.10^6$ [B-C-E-M] the right hand side is $0.16$, and when $x = v(5)$ 
a good estimate for $\log \log (x)$ is $14.10^4$. In other words, Theorem 3 
leaves far enough room for Vandiver's conjecture to be wrong!

\bigskip

\noindent {\bf References}

\bigskip

\item{[A]} D. Arlettaz: The Hurewicz homomorphism in algebraic $K$-theory, 
\break {\it J.P.A.A.} {\bf 71} (1991), 1-12.

\smallskip

\item{[As]} Ash, A.: Polyhedral reduction theory in self-adjoint
cones, in {\it Smooth Compactification of Locally Symmetric Varieties, Lie Group
s: History, Frontiers, and Applications, vol.~IV},
 by Ash, A., Mumford, D., Rapoport, M. and Tai, Y., Math. Sci. Press, Brookline,
 Mass., 1975.

\smallskip

\item{[B-K]} S. Bloch and K. Kato: $L$-functions and Tamagawa
numbers of motives. {\it In} The Grothendieck Festschrift, vol.~1,
{\it Progr. Math.} {\bf 86} (1990), Birkh\"auser, Boston, 333-400.

\smallskip

\item{[B]} A. Borel: Stable real cohomology of arithmetic groups,
{\it Ann. Scient. Ec. Norm. Sup.} 4\`eme s\'erie, t.~7
(1974), 235-272.

\smallskip

\item{[B-S]} A. Borel and J-P. Serre: Corners and arithmetic groups, {\it 
Comment. Math. Helv.} {\bf 48} (1974), 244-297.

\smallskip

\item{[Br]} K.S. Brown: Cohomology of groups, {\it Graduate Text in Math.} 
{\bf 87}, \break Springer-Verlag.

\smallskip

\item{[B-C-E-M]} J. Buhler, R. Crandall, R. Ernvall and T. Mets\"ankyl\"a: 
Irregular primes and cyclotomic invariants up to four million, {\it Math.
Comp.} {\bf 61} (1993), 151-153.

\smallskip

\item{[C-S]} J.H. Conway and N.J.A. Sloane: Sphere Packings, Lattices and 
Groups, {\it Grundlehren} {\bf 290}, Springer-Verlag.

\smallskip

\item{[D-F]} W. Dwyer and E. Friedlander: Algebraic and etale
$K$-theory, {\it Trans. Amer. Math. Soc.} {\bf 272} (1985), 247-280.

\smallskip

\item{[F]} F.T. Farrell: An extension of Tate cohomology to a class of infinite 
groups, {\it J.Pure and Applied Alg.} {\bf 10} (1977), 153-161.

\smallskip

\item{[K-L]} G.A. Kabatiansky, V.I. Levenshtein: Bounds for packings on a 
sphere and in space, {\it Problems of Information Transmission} {\bf 14} 
(1978), 1-17.

\smallskip

\item{[K]} M. Kurihara: Some remarks conjectures about cyclotomic fields and 
$K$-groups of $\Zb$, {\it Compositio Math.} {\bf 81} (1992), 223-236.

\smallskip

\item{[L]} S. Lang: Introduction to Arakelov theory, (1988),
Springer-Verlag.

\smallskip

\item{[L-S]} R. Lee and R.H. Szczarba: On the torsion in $K_4 (\Zb)$
and $K_5 (\Zb)$, {\it Duke Math. Journal} {\bf 45} No.~1
(1978), 101-130, with an Addendum by C.~Soul\'e, 131-132.

\smallskip

\item{[M]} J. Martinet: Les r\'eseaux parfaits des espaces euclidiens, (1996), 
Masson.

\smallskip

\item{[M-W]} B. Mazur and A. Wiles: Class fields of abelian
extensions of $\Qb$, {\it Invent. Math.} {\bf 76} (1984), 179-330.

\smallskip

\item{[M-H]} J. Milnor and D. Husemoller: Symmetric bilinear forms, 
{\it Ergebnisse} 
{\bf 73} (1973), Springer-Verlag.

\smallskip

\item{[S1]} C. Soul\'e: $K$-Th\'eorie des anneaux d'entiers de
corps de nombres et cohomologie \'etale, {\it Invent. Math.} {\bf
55} (1979), 251-295.

\smallskip

\item{[S2]} C. Soul\'e:
Op\'erations en K-th\'eorie alg\'ebrique,
{\it Canad. Journal of Math.} {\bf 37}, (1985), 488-550.

\smallskip

\item{[S3]} C. Soul\'e: El\'ements cyclotomiques en K-th\'eorie,
                {\it Ast\'erisque} {\bf  147-148}, (1987), 225-257.

\smallskip

\item{[S4]} C. Soul\'e, D.Abramovich, J.-F.Burnol and J.Kramer: 
{\it Lectures on Arakelov
geometry}, Cambridge Studies in Advanced Mathematics
{\bf 33},  (1992), Cambridge University Press.

\smallskip

\item{[S5]} C. Soul\'e:
On the 3-torsion in   $K_{4} (\Zb)$, {\it Topology}, to appear.

\smallskip

\item{[Su]} A.A. Suslin: Stability in algebraic $K$-theory, {\it
Lecture Notes in Math.} No.~966 (1982), 344-356.

\smallskip

\item{[V]} G. Voronoi: Nouvelles applications des param\`etres
continus \`a la th\'eorie des formes quadratiques I, {\it Crelle}
{\bf 133} (1907), 97-178.

\smallskip

\item{[W]} L.C. Washington: Introduction to cyclotomic fields, $2^{\rm nd}$ 
edition, {\it Graduate Text in Maths.} {\bf 83} (1996), Springer-Verlag.

\smallskip

\item{[Z]} S. Zhang: Positive line bundles on arithmetic surfaces, {\it Annals 
of Math.} {\bf 136} (1992), 569-587.

\bigskip

C.N.R.S. and I.H.E.S., 35 Route de Chartres, 91440, Bures-sur-Yvette, France.

\bye